\documentclass[12pt]{amsart}
\usepackage{fullpage,url,amssymb}

\DeclareFontEncoding{OT2}{}{} 

\usepackage{color}

\newcommand{\defi}[1]{\textsf{#1}} 

\newcommand{\Aff}{{\mathbb A}}

\newcommand{\F}{{\mathbb F}}

\newcommand{\PP}{{\mathbb P}}
\newcommand{\Q}{{\mathbb Q}}

\newcommand{\Z}{{\mathbb Z}}

\newcommand{\Adeles}{{\mathbf A}}

\newcommand{\calE}{{\mathcal E}}

\newcommand{\calL}{{\mathcal L}}

\newcommand{\OO}{{\mathcal O}}

\DeclareMathOperator{\ev}{ev}

\DeclareMathOperator{\Char}{char}
\DeclareMathOperator{\inv}{inv}

\DeclareMathOperator{\Hom}{Hom}

\DeclareMathOperator{\Br}{Br}

\DeclareMathOperator{\Sym}{Sym}

\DeclareMathOperator{\Pic}{Pic}

\newcommand{\injects}{\hookrightarrow}
\newcommand{\isom}{\simeq}

\newcommand{\Intersection}{\bigcap} 
\newcommand{\intersect}{\cap} 
\newcommand{\Union}{\bigcup} 
\newcommand{\union}{\cup} 
\newcommand{\tensor}{\otimes}
\newcommand{\directsum}{\oplus} 

\newtheorem{theorem}{Theorem}[section]
\newtheorem{lemma}[theorem]{Lemma}
\newtheorem{corollary}[theorem]{Corollary}

\theoremstyle{definition}

\theoremstyle{remark}
\newtheorem{remark}[theorem]{Remark}

\usepackage[
        backref,
        pdfauthor={Bjorn Poonen}, 
]{hyperref}
\usepackage[alphabetic,backrefs,lite]{amsrefs} 

\begin{document}

\title{The set of non-squares in a number field is diophantine}
\author{Bjorn Poonen}
\thanks{This research was supported by NSF grant DMS-0301280.}
\address{Department of Mathematics, University of California, 
        Berkeley, CA 94720-3840, USA}
\email{poonen@math.berkeley.edu}
\urladdr{http://math.berkeley.edu/\~{}poonen/}
\date{December 18, 2007}

\begin{abstract}
Fix a number field $k$.
We prove that $k^\times - k^{\times 2}$ is diophantine over $k$.
This is deduced from a theorem that for a nonconstant
separable polynomial $P(x) \in k[x]$, there are at most finitely many
$a \in k^\times$ modulo squares such that there is a Brauer-Manin
obstruction to the Hasse principle for the conic bundle $X$
given by $y^2 - az^2 = P(x)$.
\end{abstract}

\maketitle

\section{Introduction}\label{S:introduction}

Throughout, let $k$ be a global field;
occasionally we impose additional conditions on its characteristic.
Warning: we write $k^n = \prod_{i=1}^n k$ 
and $k^{\times n} = \{a^n:a \in k^\times\}$.

\subsection{Diophantine sets}

A subset $A \subseteq k^n$ is \defi{diophantine over $k$}
if there exists a closed subscheme $V \subseteq \Aff^{n+m}_k$
such that $A$ equals the projection of $V(k)$ under $k^{n+m} \to k^n$.
The complexity of the collection of diophantine sets over a field $k$
determines the difficulty of solving polynomial equations over $k$.
For instance, it follows from \cite{Matiyasevich1970} 
that if $\Z$ is diophantine over $\Q$,
then there is no algorithm to decide whether a
multivariable polynomial equation with rational coefficients
has a solution in rational numbers.
Moreover, diophantine sets can built up from other diophantine sets.
In particular, diophantine sets over $k$ are closed under
taking finite unions and intersections.
Therefore it is of interest to gather a library of diophantine sets.

\subsection{Main result}

Our main theorem is the following:
\begin{theorem}
\label{T:main}
For any number field $k$,
the set $k^{\times} - k^{\times 2}$ is diophantine over $k$.
\end{theorem}

In other words, there is an algebraic family 
of varieties $(V_t)_{t \in k}$ such that $V_t$ has a $k$-point
if and only if $t$ is {\em not} a square.
This result seems to be new even in the case $k=\Q$.

\begin{corollary}
For any number field $k$ and for any $n \in \Z_{\ge 0}$, the 
set $k^\times - k^{\times 2^n}$ is diophantine over $k$.
\end{corollary}

\begin{proof}
Let $A_n = k^\times - k^{\times 2^n}$.
We prove by induction on $n$ that $A_n$ is diophantine over $k$.
The base case $n=1$ is Theorem~\ref{T:main}.
The inductive step follows from $A_{n+1} = A_1 \union \{t^2 : t \in A_n\}$.
\end{proof}

\subsection{Brauer-Manin obstruction}

The main ingredient of the proof of Theorem~\ref{T:main} 
is the fact the Brauer-Manin obstruction is the only
obstruction to the Hasse principle for 
certain Ch\^atelet surfaces over number fields,
so let us begin to explain what this means.
For each place $v$ of $k$, let $k_v$ be the completion of $k$ at $v$.
Let $\Adeles$ be the ad\`ele ring of $k$.
One says that there is a Brauer-Manin obstruction
to the Hasse principle for a projective variety $X$ over $k$
if $X(\Adeles) \ne \emptyset$ but $X(\Adeles)^{\Br} = \emptyset$.

\subsection{Conic bundles and Ch\^atelet surfaces}

Let $\calE$ be a rank-$3$ vector sheaf over a base variety $B$.
A nowhere-vanishing section $s \in \Gamma(B,\Sym^2 \calE)$
defines a subscheme $X$ of $\PP\calE$ whose fibers over $B$
are (possibly degenerate) conics.
As a special case, we may take 
$(\calE,s) = (\calL_0 \directsum \calL_1 \directsum \calL_2,s_0+s_1+s_2)$
where each $\calL_i$ is a line sheaf on $B$,
and the $s_i \in \Gamma(B,\calL_i^{\tensor 2}) \subset \Gamma(B,\Sym^2 \calE)$ 
are sections that do not simultaneously vanish on $B$.

We specialize further to the case where $B=\PP^1$, 
$\calL_0 = \calL_1 = \OO$, $\calL_2 = \OO(n)$,
$s_0=1$, $s_1=-a$, and $s_2=-\tilde{P}(w,x)$
where $a \in k^\times$ and $\tilde{P}(w,x) \in \Gamma(\PP^1,\OO(2n))$
is a separable binary form of degree $2n$.
Let $P(x):=\tilde{P}(1,x) \in k[x]$, so $P(x)$ is a separable polynomial of 
degree $2n-1$ or $2n$.
We then call $X$ the conic bundle given by
\[
        y^2 - a z^2 = P(x).
\]
A \defi{Ch\^atelet surface} is a conic bundle of this type with $n=2$,
i.e., with $\deg P$ equal to $3$ or $4$.
See also~\cite{Poonen-chatelet-preprint}.

The proof of Theorem~\ref{T:main} relies on the Ch\^atelet surface case
of the following result about families of more general conic bundles:

\begin{theorem}
\label{T:finiteness}
Let $k$ be a global field of characteristic not $2$.
Let $P(x) \in k[x]$ be a nonconstant separable polynomial.
Then there are at most finitely many classes in $k^\times/k^{\times 2}$
represented by $a \in k^\times$ such that there is a Brauer-Manin
obstruction to the Hasse principle for the conic bundle $X$
given by $y^2 - az^2 = P(x)$.
\end{theorem}

\begin{remark}
\label{R:CT and Xu Fei}
Theorem~\ref{T:finiteness} is analogous to the classical fact that
for an integral indefinite ternary quadratic form $q(x,y,z)$,
the set of nonzero integers represented by $q$ over $\Z_p$ for all $p$
but not over $\Z$ fall into finitely many classes 
in $\Q^\times/\Q^{\times 2}$.
J.-L.~Colliot-Th\'el\`ene and F.~Xu explain how to interpret and
prove this fact (and its generalization to arbitrary number fields) 
in terms of the integral Brauer-Manin obstruction:
see \cite{CT-Xu2007preprint}*{\S7}, especially Proposition~7.9
and the very end of \S7.
Our proof of Theorem~\ref{T:finiteness} shares several ideas 
with the arguments there.
\end{remark}

\subsection{Definable subsets of $k_v$ and their intersections with $k$}

The proof of Theorem~\ref{T:main} requires one more ingredient, 
namely that certain subsets of $k$
defined by local conditions are diophantine over $k$.
This is the content of Theorem~\ref{T:k-definable} below,
which is proved in more generality than needed.
By a \defi{$k$-definable subset} of $k_v^n$,
we mean the subset of $k_v^n$
defined by some first-order formula in the language of fields
involving only constants from $k$, even though the variables
range over elements of $k_v$.

\begin{theorem}
\label{T:k-definable}
Let $k$ be a number field.
Let $k_v$ be a nonarchimedean completion of $k$.
For any $k$-definable subset $A$ of $k_v^n$,
the intersection $A \intersect k^n$ is diophantine over $k$.
\end{theorem}

\subsection{Outline of paper}

Section~\ref{S:definable} shows that Theorem~\ref{T:k-definable}
is an easy consequence of known results, namely the description of 
definable subsets over $k_v$, and the diophantineness of
the valuation subring $\OO$ of $k$ defined by $v$.
Section~\ref{S:conic bundles} proves Theorem~\ref{T:finiteness}
by showing that for most twists of a given conic bundle,
the local Brauer evaluation map at one place is enough to rule out
a Brauer-Manin obstruction.
Finally, Section~\ref{S:diophantine}
puts everything together to prove Theorem~\ref{T:main}.
\
\section{Subsets of global fields defined by local conditions}
\label{S:definable}

\begin{lemma}
\label{L:m-th powers}
Let $m \in \Z_{>0}$ be such that $\Char k \nmid m$.
Then $k_v^{\times m} \intersect k$ is diophantine over $k$.
\end{lemma}

\begin{proof}
The valuation subring $\OO$ of $k$ defined by $v$ is diophantine over $k$: 
see the first few paragraphs of \S3 of \cite{Rumely1980}.
The hypothesis $\Char k \nmid m$ implies
the existence of $c \in k^\times$ 
such that $1+c\OO \subset k_v^{\times m}$; fix such a $c$.
The denseness of $k^\times$ in $k_v^\times$ implies
$k_v^{\times m} \intersect k = (1+c\OO) k^{\times m}$.
The latter is diophantine over $k$.
\end{proof}

\begin{proof}[Proof of Theorem~\ref{T:k-definable}]
Call a subset of $k_v^n$ \defi{simple}
if it is of one of the following two types: 
$\{\vec{x} \in k_v^n: f(\vec{x})=0\}$
or $\{\vec{x} \in k_v^n: f(\vec{x}) \in k_v^{\times m}\}$
for some $f \in k[x_1,\ldots,x_n]$ and $m \in \Z_{>0}$.
It follows from the proof of \cite{Macintyre1976}*{Theorem~1}
(see also \cite{Macintyre1976}*{\S2} and \cite{Denef1984}*{\S2})
that any $k$-definable subset $A$ is a boolean combination of simple subsets.
The complement of a simple set of the first type is a simple set of
the second type (with $m=1$).
The complement of a simple set of the second type is a union of simple sets,
since $k_v^{\times m}$ has finite index in $k_v^\times$.
Therefore any $k$-definable $A$ is a finite union of finite intersections
of simple sets.
Diophantine sets in $k$ are closed under taking finite unions 
and finite intersections,
so it remains to show that 
for every simple subset $A$ of $k_v^n$,
the intersection $A \intersect k$ is diophantine.
If $A$ is of the first type, then this is trivial.
If $A$ is of the second type, then this follows 
from Lemma~\ref{L:m-th powers}.
\end{proof}

\section{Family of conic bundles}
\label{S:conic bundles}

For a place $v$ of $k$
let $\Hom'(\Br X,\Br k_v)$ be the set of $f \in \Hom(\Br X,\Br k_v)$
such that the composition $\Br k \to \Br X \stackrel{f}\to \Br k_v$ 
equals the map induced by the inclusion $k \injects k_v$.
The $v$-adic evaluation pairing $\Br X \times X(k_v) \to \Br k_v$
induces a map $X(k_v) \to \Hom'(\Br X,\Br k_v)$.

\begin{lemma}
\label{L:surjective}
With notation as in Theorem~\ref{T:finiteness},
there exists a finite set of places $S$ of $k$,
depending on $P(x)$ but not $a$,
such that if $v \notin S$ and $v(a)$ is odd,
then $X(k_v) \to \Hom'(\Br X,\Br k_v)$ is surjective.
\end{lemma}

\begin{proof}
The function field of $\PP^1$ is $k(x)$.
Let $k(X)$ be the function field of $X$.
Let $Z$ be the zero locus of $\tilde{P}(w,x)$ in $\PP^1$.
Let $G$ be the group of $f \in k(X)^\times$
having even valuation at every closed point of $\PP^1-Z$.
Choose $P_1(x),\ldots,P_m(x) \in G$ representing a $\F_2$-basis 
for the image of $G$ in $k(x)^\times/k(x)^{\times 2} k^\times$.
We may assume that $P_m(x)=P(x)$.
Choose $S$ so that each $P_i(x)$ is a ratio of polynomials 
whose nonzero coefficients are $S$-units.
A well-known calculation (see \cite{Skorobogatov2001}*{\S7.1}) shows
that the class of each quaternion algebra $(a,P_i(x))$ in $\Br k(X)$
belongs to the subgroup $\Br X$,
and that the cokernel of $\Br k \to \Br X$
is an $\F_2$-vector space with the classes of $(a,P_i(x))$ for $i \le m-1$
as a basis.

Suppose that $v \notin S$ and $f \in \Hom'(\Br X,\Br k_v)$.
The homomorphism $f$ is determined by 
where it sends $(a,P_i(x))$ for $i \le m-1$.
We need to find $R \in X(k_v)$ mapping to $f$.

Let $\OO_v$ be the valuation ring in $k_v$, 
and let $\F_v$ be its residue field.
We may assume that $\Char \F_v \ne 2$.
For $i \le m-1$, 
choose $c_i \in \OO_v^\times$ whose image in $\F_v^\times$ is a square or not,
according to whether $f$ sends $(a,P_i(x))$ to $0$ or $1/2$
in $\Q/\Z \isom \Br k_v$.
Since $v(a)$ is odd, we have $(a,c_i)=(a,P_i(x))$ in $\Br k_v$.

View $\PP^1-Z$ as a smooth $\OO_v$-scheme,
and $Y$ be the finite \'etale cover of $\PP^1-Z$
whose function field is obtained by adjoining $\sqrt{c_i P_i(x)}$ 
for $i \le m-1$ and also $\sqrt{P(x)}$.
Then the generic fiber $Y_{k_v} := Y \times_{\OO_v} k_v$ 
is geometrically integral.
Assuming that $S$ was chosen to include all $v$ with small $\F_v$,
we may assume that $v \notin S$ implies that $Y$ has a (smooth) $\F_v$-point,
which by Hensel's lemma lifts to an $k_v$-point $Q$.
There is a morphism from $Y_{k_v}$
to the smooth projective model of $y^2=P(x)$ over $k_v$,
which in turn embeds as a closed subscheme of $X_{k_v}$, 
as the locus where $z=0$.
Let $R$ be the image of $Q$ under $Y(k_v) \to X(k_v)$, 
and let $\alpha=x(R) \in k_v$.
Evaluating $(a,P_i(x))$ on $R$ yields $(a,P_i(\alpha))$,
which is isomorphic to $(a,c_i)$
since $c_i P_i(\alpha) \in k_v^{\times 2}$.
Thus $R$ maps to $f$, as required.
\end{proof}

\begin{proof}[Proof of Theorem~\ref{T:finiteness}]
Let $S$ be as in Lemma~\ref{L:surjective}.
Enlarge $S$ to assume that $\Pic \OO_{k,S}$ is trivial.
Then the set of $a \in k^\times$ such that $v(a)$ is even for all $v \notin S$
has the same image in $k^\times/k^{\times 2}$ as 
the finitely generated group $\OO_{k,S}^\times$,
so the image is finite.
Therefore it will suffice to show that if $v \notin S$ and $v(a)$ is odd,
then the corresponding surface $X$ has no Brauer-Manin obstruction
to the Hasse principle.

If $X(\Adeles)=\emptyset$, then the Hasse principle holds.
Otherwise pick $Q = (Q_w) \in X(\Adeles)$, where $Q_w \in X(k_w)$ for each $w$.
For $A \in \Br X$, let $\ev_A\colon X(L) \to \Br L$ be the evaluation map
for any field extension $L$ of $k$, and let $\inv_w\colon \Br k_w \to \Q/\Z$
be the usual inclusion map.
Define
\begin{align*}
        \eta\colon \Br X &\to \Q/\Z \isom \Br k_v \\
        A &\mapsto - \sum_{w \ne v} \inv_w \ev_A(Q_w).
\end{align*}
By reciprocity, $\eta \in \Hom'(\Br X,\Br k_v)$.
By Lemma~\ref{L:surjective}, there exists $R \in X(k_v)$ giving rise to $\eta$.
Define $Q' = (Q'_w) \in X(\Adeles)$ by $Q'_w:=Q_w$ for $w \ne v$ and $Q'_v:=R$.
Then $Q' \in X(\Adeles)^{\Br}$, so there is no Brauer-Manin obstruction
to the Hasse principle for $X$.
\end{proof}

\section{The set of nonsquares is diophantine}
\label{S:diophantine}

\begin{proof}[Proof of Theorem~\ref{T:main}]
For each place $v$ of $k$, 
define $S_v := k^\times \intersect k_v^{\times 2}$
and $N_v:=k^\times - S_v$.
By Theorem~\ref{T:k-definable}, 
the sets $S_v$ and $N_v$ are diophantine over $k$.

By \cite{Poonen-chatelet-preprint}*{Proposition~4.1},
there is a Ch\^atelet surface
\[
        X_1 \colon y^2 - b z^2 = P(x)
\]
over $k$, with $P(x)$ a product of two irreducible quadratic polynomials,
such that there is a Brauer-Manin obstruction to the Hasse principle for $X_1$.
For $t \in k^\times$, let $X_t$ be the (smooth projective) Ch\^atelet surface 
associated to the affine surface
\[
        U_t\colon y^2 - tb z^2 = P(x).
\]

We claim that the following are equivalent for $t \in k^{\times}$:
\begin{enumerate}
\item[(i)] $U_t$ has a $k$-point.
\item[(ii)] $X_t$ has a $k$-point.
\item[(iii)] $X_t$ has a $k_v$-point for every $v$ 
and there is no Brauer-Manin obstruction to the Hasse principle for $X_t$.
\end{enumerate}
The implications (i)$\implies$(ii)$\implies$(iii) are trivial.
The implication (iii)$\implies$(ii) 
follows from \cite{CT-Coray-Sansuc1980}*{Theorem~B}.
Finally, 
in \cite{CT-Coray-Sansuc1980}, the reduction of Theorem~B to Theorem~A
combined with Remarque~7.4 shows that (ii) 
implies that $X_t$ is $k$-unirational, which implies~(i).

Let $A$ be the (diophantine) set of $t \in k^{\times}$ such that (i) holds.
The isomorphism type of $U_t$ 
depends only on the image of $t$ in $k^\times/k^{\times 2}$,
so $A$ is a union of cosets of $k^{\times 2}$ in $k^\times$.
We will compute $A$ by using~(iii).

The affine curve $y^2 = P(x)$ is geometrically integral
so it has a $k_v$-point for all places $v$ outside a finite set $F$.
So for any $t \in k^\times$,
the variety $X_t$ has a $k_v$-point for all $v \notin F$.
Since $X_1$ has a $k_v$-point for all $v$ and in particular for $v \in F$, 
if $t \in \Intersection_{v \in F} S_v$,
then $X_t$ has a $k_v$-point for all $v$.

Let $B:=A \union \Union_{v \in F} N_v$.
If $t \in k^\times - B$, then $X_t$ has a $k_v$-point for all $v$,
and there is a Brauer-Manin obstruction to the Hasse principle for $X_t$.
By Theorem~\ref{T:finiteness}, $k^\times-B$ consists of finitely many
cosets of $k^{\times 2}$, one of which is $k^{\times 2}$ itself.
Each coset of $k^{\times 2}$ is diophantine over $k$,
so taking the union of $B$ with all the finitely many missing cosets
except $k^{\times 2}$ shows that $k^\times - k^{\times 2}$
is diophantine.
\end{proof}


\section*{Acknowledgements} 

I thank Jean-Louis Colliot-Th\'el\`ene for a few comments,
and Alexandra Shlapentokh for suggesting some references.

\end{document}